\documentclass[twoside,a4paper]{article}

\usepackage{amsmath,graphicx,amssymb,fancyhdr,amsthm,enumerate,textcomp}
\usepackage[usenames]{color}
\newtheorem{thm}{Theorem}[section]

\newtheorem{prop}[thm]{Proposition} 
 
\theoremstyle{definition}

\theoremstyle{remark} 

\def\beq{\begin{eqnarray}} 
\def\eeq{\end{eqnarray}} 
\def\bsp{\begin{split}} 
\def\esp{\end{split}}

\def\d{\mathrm{d}}

\def \epsilont {\tilde{\epsilon}} \def \betat {\tilde{\beta}} \def \alphat {\tilde{\alpha}} \def \gammat
{\tilde{\gamma}} \def \kappat {\tilde{\kappa}} \def \sigmat {\tilde{\sigma}} \def \rhot {\tilde{\rho}} \def \taut {\tilde{\tau}}

\begin{document}   
   
\title{\Large\textbf{Mathematical Properties of a Class of Four-dimensional Neutral Signature Metrics}}
\author{{\large\textbf{D. Brooks$^{1}$,~N. Musoke$^{2}$, ~ D. McNutt$^{1}$,~ and~ A. Coley$^{1}$ }} \vspace{0.3cm} \\
$^{1}$Department of Mathematics and Statistics,\\ Dalhousie University, Halifax, Nova Scotia,\\ Canada B3H 3J5 \\ 
$^{2}$ Perimeter Institute For Theoretical Physics, \\ Waterloo, Ontario \\
 \texttt{dario.a.brooks,ddmcnutt,msknathan@dal.ca, aac@mathstat.dal.ca} }

\date{\today} \maketitle \pagestyle{fancy} \fancyhead{} \fancyhead[CE]{Brooks et al.} \fancyhead[LE,RO]{\thepage}
\fancyhead[CO]{Neutral signature metrics} \fancyfoot{} 

\begin{abstract} While the Lorenzian and Riemanian metrics for which all 
polynomial scalar curvature invariants vanish (the VSI property) are well-studied, less is known about the 
four-dimensional neutral signature metrics with the VSI property. Recently it was shown that the neutral signature metrics belong to two distinct 
subclasses: the Walker and Kundt metrics. In this paper we have chosen an example from each of the two subcases of the Ricci-flat VSI Walker metrics respectively. 

To investigate the difference between the metrics we determine the existence of a null, geodesic, expansion-free, shear-free and vorticity-free vector,
and classify these spaces using their holonomy algebras. The geometric implications of these algebras are further studied by 
identifying the recurrent or covariantly constant null vectors, whose existence is required by the holonomy structure in each example. We conclude 
the paper with a simple example of the equivalence algorithm for these pseudo-Riemannian manifolds, which is the only approach to classification that provides all necessary information to determine equivalence. \vspace{0.3mm}

 \end{abstract}
   
\section{Introduction}

Let us consider the four-dimensional pseudo-Riemannian spaces of signature $(2,2)$, the so called neutral signature; in this paper we will study the mathematical properties of the neutral signature solutions for which all of polynomial invariants formed by the curvature tensor and its covariant derivatives vanish. If all polynomial scalar curvature invariants vanish, we say the space has the {\it VSI property} and is hence a {\it VSI space} \cite{Alcolado};  in analogy with the 4D degenerate Kundt metrics with the VSI property \cite{degen, VSI0} in the Lorentzian case. 

Through the investigation of how the Riemann tensor and its covariant derivatives change under boosts in each of the null planes, it was possible to classify the pseudo-Riemannian spaces through the boost-weight decomposition \cite{SIGII, SIGIII}. It was shown that if a frame is found where the curvature tensors of the space have only negative boost-weight terms (the {\bf N} property) then it is a VSI spacetime. Furthermore, from this result it was shown that the VSI spacetimes were either Kundt (possessing a geodesic, expansion-free, shear-free, and twist-free null-congruence) or Walker (admiting a 2-dimensional invariant null plane \cite{Walker}) form \cite{SIGIII}. In \cite{CCMMB} we presented examples of  4D neutral signature VSI metrics which are genuinely Walker spaces (i.e., not Kundt).  

To illustrate this dichotomy in the neutral Ricci-flat VSI Walker metrics, we study two distinct subcases, one which is Walker in general \cite{CCMMB,Alcolado}, 
and another that is strictly Kundt. We compare these metrics by determining the existence of a null geodesic, expansion-free, shear-free, and vorticity-free 
vector using the spin-coefficient formalism \cite{Law}, and then use the Lie algebra classification provided in \cite{GT} and \cite{Wang} to distinguish the two metrics. We show that this classification is well-suited for determining the existence of covariant constant null vectors, recurrent null vectors and more general invariant null distributions, although these comparisons are only helpful for showing when two metrics are not equivalent.

A natural question is to ask when, if at all, are the subcases of the two metrics equivalent. This question can be answered by
implementing the equivalence algorithm for neutral metrics. We end this paper with an example of the equivalence algorithm
applied to a simple subcase of the Kundt-Walker metric. 

\section{4D Neutral signature Walker metrics}

A metric is said to posses a 2-dimensional invariant plane if there exist null vectors $l$ and $m$ such that the bivector
$l\wedge m$ is recurrent; that is, \begin{equation} \nabla_a (l\wedge m ) = k_a (l\wedge m) \end{equation} for some covector
$k_a$. If these vectors are also null and orthogonal, the invariant plane is called totally null. 4D neutral signature spaces
possessing an invariant null plane are known as Walker metrics \cite{Walker}. In \cite{Law}, it was shown that a metric has an
invariant null plane if and only if there exists a frame in which the spin coefficients $\kappa = \rho = \sigma=\tau=0$.

In analogy with the Lorentzian case, we will say a metric is Kundt if it possesses a non-zero null vector $\ell$ which is geodesic, expansion-free, twist-free, and
shear-free, which implies a particular form for the covariant derivative of $\ell$ \cite{SIGIII}. This condition implies that
there exists a null coframe, $\{ n_a, \ell_a, \bar{m}_a, m_a\}$, in which the spin-coefficients are required to satisfy:
 \beq \kappat =
\kappa = \rhot = \rho = \sigmat =\sigma = 0. \nonumber \eeq \noindent Equivalently, if this is the case, the covariant derivative of $\ell_a$ is of the form
\beq \ell_{a;b} = -(\epsilon' - \epsilont') \ell_a \ell_b + (\alphat' + \beta') \ell_a m_b - \taut m_a \ell_b \nonumber \eeq
\noindent which is automatically geodesic,  expansion-free, shear-free and
vorticity-free. 

We are especially interested in the case where the VSI metric is Walker, and hence admits an invariant null plane, but is
not Kundt. For Walker metrics in 4D neutral spaces with an invariant 2-dimensional null plane it is always possible to find a
null field that is geodesic, expansion-free, vorticity-free (see below). In general, this null vector is not necessarily
shear-free.  To achieve this, we work with a particular class of VSI-Walker Ricci-flat metrics \cite{CCMMB} in the Walker form: \beq \mathrm{d}s^2=2\d u(\d v+A\d u+C\d
U)+2\d U(\d V+B\d U), \label{Walker2}\eeq where \beq A&=&vA_1(u,U)+VA_2(u,U)+A_0(u,U), \nonumber \\ B&=&VB_{10}(u,U) +
v^2B_{02}(u,U)+vB_{01}(u,U)+B_{00}(u,U) \nonumber \\ C&=&vC_{11}(u,U)+VC_2(u,U)+C_0(u,U) \nonumber \eeq 
\noindent and the remaining arbitrary functions in the metric are assumed to be analytic.

For the space to be Ricci flat, we must also have that $A_2 B_{02}=0$. This metric does not in general possess the ${\bf N}$-property, but rather
the weaker requirement of the ${\bf N}^G$-property: where the boost-weights of the Riemann tensor may be be treated as vectors and transformed into boost-weight vectors satisfying the {\bf N} property \cite{SIGII, SIGIII}.

The null tetrad frame $\{\ell, n, m, -\tilde{m}\}=\{l_1, n_1, l_2, n_2\}$ is defined by: \beq l_1 &= du, ~~n_1= dv+A(u,v,U,V) du
+ \frac{C(u,v,U,V)}{2}dU \label{Coframe1} \\ l_2 &= dU, ~~n_2 = dV +\frac{C(u,v,U,V)}{2}du+ B(u,v,U,V)dU \label{Coframe2} \eeq
The invariant null plane is given by the null orthogonal vectors $l_1$ and $l_2$: \begin{equation} \nabla \left( l_1 \wedge l_2
\right) = \left( l_1 \wedge l_2 \right) \otimes \left( \left( \frac{C_2}{2}+A_1 \right) du + \left( B_{10} + \frac{C_{11}}{2}
\right) dU \right) \end{equation} confirms that this is a Walker metric.

There are two subcases for $A_2 B_{02}=0$. The Walker-Kundt case {$B_{02} \neq 0; A_2=0$} was investigated in \cite{CCMMB}.
Ricci flat solutions in the case $A_2= C_2=0$, $B_{02} \neq 0$ were obtained, and the corresponding holonomy
algebra was found to be $A_{26}$ \cite{GT}, and so our metric has a null two-dimensional distribution containing a recurrent
vector field. When $B_{10,u}=0$, the holonomy algebra reduces to $A_{17}$, and so the metric has a null
two-dimensional distribution containing one parallel vector field. 


\subsection{Kundt Condition for Walker VSI Metrics}

Despite the difference in the geometric structure of the neutral metrics, the definitions of twist, expansion and shear may be generalized to the neutral case. While a physical interpretation analogous to the work of Ehlers, Sachs and Kundt in the Lorentzian case is no longer available, the computation of these quantities still relies on projecting them onto a hyper-surface orthogonal to the null direction. However, this hyper-surface will now be timelike. 

As an example, in the case that these quantities vanish for a particular null direction, the surface orthogonal to this null direction ($u=constant$, or $U=constant$) may be treated as a two-dimensional Lorentzian manifold. This is a special property, which does not occur for all Ricci-flat $VSI$ Walker metrics, due to non-vanishing shear. To illustrate this point, we explicitly show that the $VSI$ Walker metric given in Theorem 2.1 of \cite{CCMMB} generally do not admit a null vector $X_a$  which is geodesic, expansion-free ($X^{a}_{;a} = 0$),  shear-free ($X^{a;b}X_{(a;b)}=0$) and vorticity-free ($X^{a,b}X_{[a,b]}=0$).

\begin{prop} \label{WalkerKundt} A four dimensional neutral signature Walker metric of the form: 
\beq 
\mathrm{d}s^2=2\d u(\d v+A\d u+C\d U)+2\d U(\d V+B\d U), 
\label{Walker20}\eeq 
\noindent with $A,B,C$ of the form:
\beq
A&=&vA_1(u,U)+VA_2(u,U)+A_0(u,U), \nonumber \\
B&=&VB_1(u,v,U)+B_0(u,v,U)  \nonumber \\
C&=&C_1(u,v,U)+VC_2(u,U)+C_0(u,U), \nonumber 
\eeq
\noindent  is Kundt with a null geodesic, expansion-free, shear-free and vorticity-free vector $X_a$ which is proportional to: \begin{itemize} \item $\ell_a$ if
and only if the metric component $A$ in $\ell_1$ in \eqref{Coframe1} satisfies $A_{,V} = 0$. \item $m_a$ if and only if the metric component $B$ in $\ell_2$ in \eqref{Coframe2} satisfies $B_{,v} = B_{,v~v} = 0$. \end{itemize} \end{prop}

\begin{proof}

As in the Lorentzian signature case, we find that that the null geodesic vector must satisfy, 
\beq \frac12 \epsilon^{abcd}X_b X_{c;d} = \omega X^a, \nonumber \eeq

\noindent where $\omega^2 = X_{[a;b]} X^{a;b}$. Imposing the vorticity-free condition, the above implies $X_a$ is hypersurface orthogonal. Locally, we may choose this vector field to be the gradient of some function $X(u,v,U,V)$. Then by expanding and equating powers of $v$ and $V$ we find that the conditions $X$ be null and geodesic both imply that $X_{,v} = X_{,V} = 0$. Thus $X_{,a} = X_{,U} \ell_a + X_{,u} m_a $ and its covariant derivative is of the form:
\small \beq X_{a;b} &=& (A_2 X_{,U} + A_1 X_{,u} + X_{,uu}) m_a m_b + (C_2 X_{,u} + C_{11} X_{,u} + X_{,Uu}) [\ell_a m_b + m_a
\ell_b] \nonumber \\ && + (X_{,U} B_{10} + 3 X_{,u} v^2 B_{03} + 2 X_{,u} v B_{02} + X_{,u} B_{01} + X_{,UU}) \ell_a \ell_b.
\nonumber \eeq
\normalsize 

Given the vector field $X^a$ we may complete the basis for the tangent space, $\{ X^a, Y^a, m^{a},\tilde{m}^{a}\}$ with $X_aY^a = 1$, and assuming $X^a$ is null, geodesic, expansion-free, shear-free, and vorticity free we may write the covariant derivative of $X^a$ as 

\beq X_{a;b} = L_{11} X_a X_b + L_{12} X_a m_b + \tilde{L}_{12} X_a \tilde{m}_b. \label{KundtVect} \eeq 

\noindent To identify the terms that must vanish in the covariant derivative of $X^a$ relative to the original basis for the tangent space, we  project this tensor to the subspace of the tangent space perpendicular to $X^a$. For the moment, we will assume that $X_{,U} \neq 0$ and consider the following projection operator:

\beq h_a^{~b} = g_a^{~b} + \frac{\tilde{m}_a X^b}{X_{,U}} + \frac{X_a \tilde{m}^b}{X_{,U}}. \nonumber \eeq
\noindent As $h_a^{~b} X_b = 0$ this will serve to project onto the subspace $\{e_{i'} \}$. Applying this operator to the
covariant derivative of $X_{a}$, $X_{c;d}h_a^{~c} h_b^{~d}$, we have one non-zero coefficient of $m_a m_b$, namely the component
$X_{c;d}h_a^{~c} h_b^{~d}\tilde{m}^a \tilde{m}^b$ which is of the form: \small

\beq & X^{-2}_{,U} [X_{,u} (-2 X_{,u} X_{,U} V B_{11} - X_{,u} X_{,U} B_{10} - 3X_{,u}^2 V^2 B_{03} - 2 X_{,u}^2 V B_{02} ] &
\nonumber \\ & +X^{-2}_{,U}[- X_{,u}^2 B_{01} - X_{,u} X_{,UU} + 2X_{,U}^2 C_2 + 2X_{,u} X_{,U}C_{11} ] & \nonumber \\ & +
X^{-2}_{,U}[ 2 X_{,U} X_{,Uu} - A_1 X_{,U}^2]- A_2 X_{,U}^3 - X_{,uu} X_{,U}^2 & \nonumber \eeq

\normalsize

\noindent Equating the $v$ and $V$ linear terms we must have $X_{,u} = 0$ or $B_{11} = B_{03} = B_{02} =0$. For the moment we
assume that $X_{,u} =0$, producing the simpler expression for this coefficient of $m_a m_b$, the vanishing of which requires that $A_2$ must vanish and $X_a$ be proportional to $\ell_a$ 


Alternatively, assuming $X_{,u} \neq 0$ we may produce a similar projection operator, and then by contracting this new
projection operator twice with the covariant derivative of $X_{a;b}$ we find a similar expression to that above 

\beq & X^{-2}_{,u} [X_{,u} (-2 X_{,u} X_{,U} V B_{11} - X_{,u} X_{,U} B_{10} - 3X_{,u}^2 V^2 B_{03} - 2 X_{,u}^2 V B_{02} ] &
\nonumber \\ & +X^{-2}_{,U}[- X_{,u}^2 B_{01} - X_{,u} X_{,UU} + 2X_{,U}^2 C_2 + 2X_{,u} X_{,U}C_{11} ] & \nonumber \\ & +
X^{-2}_{,U}[ 2 X_{,U} X_{,Uu} - A_1 X_{,U}^2]- A_2 X_{,U}^3 - X_{,uu} X_{,U}^2 & \nonumber \eeq

\noindent As before, the vanishing of the $v$ and $V$ linear terms imply that either $X_{,U } =0$ or $B_{11}=B_{03} = B_{02} =0$.
It is not clear if there is a solution to the partial differential equation for $X(u,U)$ when these three functions vanish, thus
we will assume that $X_{,U} =0$. Simplifying the equations, we find that this metric will be Kundt if and only if $B_{01} = B_{02} =B_{03} =0$ with $X_a$ proportional to $m_a$. 
\end{proof}

\begin{subsection}{Example 1: A Non-Kundt Walker Metric}

We examine the case in which $B_{02}=0$ and $A_2$ may or may not be zero. From section 3.1, when $A_2=0$ or $B_{01}=B_{02}=0$ this is 
a Kundt metric, and so we will assume that these are non-zero in general.

To satisfy the condition of Ricci flatness, the metric functions must be solutions to the folllowing equations: 
\beq
&\label{eqn:RicciFlat1}\frac{\partial C_2}{\partial u}+\mathcal{A}C_2 = 2 \frac{\partial A_2}{\partial U}-2 \mathcal{B} A_2, &
\\ &\label{eqn:RicciFlat2}2\frac{\partial B_{01}}{\partial u}+\mathcal{A}B_{01} = \frac{\partial C_{11}}{\partial U} +
\mathcal{B} C_{11},& \\ &\frac{\partial \mathcal{A}}{\partial U} + \frac{\partial \mathcal{B}}{\partial u} - 2 A_2 B_{01} +
\frac{1}{2}C_2 C_{11} =0, & \eeq 
\noindent where $\mathcal{A}\equiv A_1-\frac{1}{2}C_2$ and $\mathcal{B}\equiv
B_{10}-\frac{1}{2}C_{11}$. As a simple example, we set $\mathcal{A} = \mathcal{B} =0$ and obtain a simple solution for
\eqref{eqn:RicciFlat1} and \eqref{eqn:RicciFlat2}: \beq &A_1 = \frac{1}{2}C_2,~~A_2 = a U + \frac{\alpha}{u},~~ B_{10} =
\frac{1}{2} C_{11},~~B_{01} = \frac{1}{2}\left(\frac{c_2}{U} + du \right),& \nonumber \\ &C_2 = 2a u + \frac{\beta}{U},~~C_{11}
= \frac{c_1}{u} + d U, & \nonumber \eeq
 \noindent where $a, \alpha, \beta, c_1, c_2, d$ are constants. Then the Ricci tensor has
one nonzero component: \begin{equation}\label{RicciB02eq0} R_{13} = R_{31} = \frac{1}{2}\frac{(\beta c_1 - 2\alpha c_2) + (2a
c_1 + \beta d - 2 c_2 a - 2 \alpha d) u U}{uU}, \end{equation} which for Ricci flatness gives us algebraic constraints, $\beta
c_1 = 2\alpha c_2$ and $2a c_1 + \beta d - 2 c_2 a - 2 \alpha d = 0$. If $c_1 \neq 0$, we must have $\beta=2 c_2 \alpha/c_1$ and
$(c_1 - c_2)(a- \alpha d/c_1)=0$. We chose to work with $a = \alpha d/c_1$, and assume $c_1$ is nonzero and $c_1 \neq c_2$ in
order to make later calculations more manageable.

When $c_1$ vanishes, the Ricci-flat conditions produce five possible subcases where some of the constants must vanish or satisfy
an identity:

\begin{itemize} \item $\alpha$ = 0, and $\beta$ = $\frac{2ac_2\alpha}{d}$. \item $\alpha$ = 0, d = 0, and a = 0. \item $\alpha$ = 0, d
= 0, and $c_2$ = 0. \item $c_2$ = 0, and d = 0. \item $c_2$ = 0, and $\beta$ = 2$\alpha$. \end{itemize}

\noindent while if $c_1$ is non-zero and $c_2 = c_1$ we find one more case:
\begin{itemize} \item $\beta = 2 \alpha$ \end{itemize} \noindent The exact form of the two recurrent vectors for these subcases
are not discussed in this paper, however, the analysis will be similar to the case studied in this paper.

\end{subsection} 

\begin{subsection}{Example 2: A Kundt Walker Metric } To provide a simple example of the equivalence
algorithm, we consider a case that is automatically Kundt and Ricci-flat, with $B_{02} \neq 0$ and $A_2 = C_2 = 0$; ensuring
that this is indeed a Walker-Kundt metric, with only one null, geodesic, expansion-free, shear-free and vorticity free vector.
Imposing the Ricci-flat conditions we have the following expressions for the metric functions: \small \beq A_{0} &=& \frac18
\frac{-2B_{10} C_{11} - 4 A_1 B_{01} +4 B_{01,u} - 2 C_{11,U} + C_{11}^2 }{B_{02}} \nonumber \\ A_{1} &=& \frac12
[log(B_{02})]_{,u} \nonumber \\ C_{11} &=& 2 B_{10} + [log B_{02}]_{,U} + G(U) \nonumber \eeq \normalsize \noindent where $A_0$
has not been fully expanded in order to display it compactly. With these metric functions, the components of the Riemann tensor
are now: \small \beq R_{1224} &=& R_{2434}= B_{10,u} \nonumber \\ R_{2323} &=& 2 B_{02} \nonumber \\ R_{2424} &=& -C_{0,U
u}+B_{00,u u}-3 v^2 A_{1} B_{02,u}-2 v^2 A_{1,u} B_{02}-B_{10} v C_{11,u}+B_{10} v A_{1,U} \nonumber \\ &&+A_{1} V
B_{10,u}-A_{1} B_{10} C_{0}+2 v^2 A_{1}^2 B_{02}+B_{10} C_{11} A_{0}+(1/2) v C_{11} C_{11,u} \nonumber \\ &&-4 A_{0} v
B_{02,u}-2 A_{0,u} v B_{02}-v A_{1} B_{01,u}-v A_{1,u} B_{01}+4 v A_{1} B_{02} A_{0}+A_{0,U U} \nonumber \\ &&+v^2 B_{02,u u}+V
B_{10,u u}-2 A_{0} B_{01,u}-A_{0,u} B_{01}-B_{10} C_{0,u}-B_{10,u} C_{0} \nonumber\\ && -A_{1} C_{0,U}-A_{1,U} C_{0}+v A_{1,U,
U}+C_{11} A_{0,U}+C_{11,U} A_{0}+2 A_{0}^2 B_{02}\nonumber \\ &&+(1/2) C_{11,u} C_{0} -v C_{11,U u}+v B_{01,u u}+A_{1}
B_{00,u}+B_{10} A_{0,U}-B_{10,u} v C_{11} \nonumber \eeq 
\normalsize
\noindent For the remainder of the paper, the component $R_{2424}$ will
be denoted as $\Psi$. It should be remarked that when the Ricci-flat conditions are imposed $\Psi$ is independent of $v$ and
$V$. In the current form it is not entirely clear that this is the case. 

\end{subsection}

\begin{section}{Holonomy Algebras for Neutral Walker Metrics}
From theorem 8.5 in \cite{Hall}, the holonomy group must preserve the inner product of the neutral metric, and so this group
will be a subgroup of the generalized orthogonal group O(2,2), and the holonomy algebras will be subalgebras of
$\mathfrak{o}(2,2)$, with members of $\phi$ represented as 2-forms. Alternatively we may represent the elements of 
$\mathfrak{o}(2,2)$ as (1,1) tensors by raising one index of the 2-form representation of each element.
 
In \cite{GT} a classification of all possible Lie subalgebras of $\mathfrak{o}(2,2)$ is given by exploiting the 
isomorphism between $\mathfrak{o}(2,2)$ and $\mathfrak{su}(1,1) \times \mathfrak{su}(1,1)$, producing 32 possible classes of Lie 
subalgebras. Using this classification, the authors were able to examine 31 of the 32 cases and determine whether the subalgebra in 
each case is achieved for a particular neutral metric in four dimensions as a holonomy algebra. The remaining case, $A_{13}$, was 
shown to be the holonomy algebra of a neutral metric \cite{Krantz,Davidov}.

The geometric structure was determined for each subalgebra and summarized in Table II of \cite{GT}. This approach relied upon the 
existence of a Lie algebra isomorphism and the known structure of $\mathfrak{su}(1,1)$ to enumerate all possibilities. In the 
neutral 4D case, if the distribution is 2D null, it need not contain an invariant null line implying that a Walker space is not
necessarily a Kundt space. If a space admits an one dimensional holonomy invariant distribution (an invariant null-line) it 
is automatically a special Kundt spacetime admitting a recurrent vector, and admits an invariant null plane containing the null 
line. Thus, a Walker space admitting an invariant one-dimensional distribution must contain an invariant 2D distribution.

There is an alternative approach based on geometric and algebraic considerations for the neutral metric manifolds. This
formalism was used in Wang and Hall \cite{Wang} to classify the holonomy subalgebras in order to study the problem of
projectively related manifolds sharing similar holonomy groups.

For an arbitrary {\it orthonormal basis} for the tangent space of M at m, $T_m M$, satsifying $g_m(x,x)=g_m(y,y) = -g_m(s,s)
=-g_m(t,t) = 1$, the elements of the six-dimensional vector space of 2-forms at m, $\Lambda_m M$ may be represented as tensors
of type (2,0), (1,1) or (0,2). Raising the indices we call members of $\Lambda_m M$ {\it bivectors} and express these in
component form as: $F \in \Lambda_m M$, $F \leftrightarrow F^{ab} (= - F^{ba})$. The bivector representation of
$\mathfrak{o}(2,2)$ is the Lie algebra $\{\alpha \in M_4 \mathbb{R} : \alpha \epsilon + (\alpha \epsilon)^T \}$ with $\epsilon =
diag(1,1,-1,-1)$ and T denoting matrix transpose. There is a natural metric $P$ on $\Lambda_m M$ for which the inner product
$P(F,G)$ of $F, G \in \Lambda_m M$ is $F^{ab}G_{ab}= P_{abcd}F^{ab}G^{cd}$, with $P_{abcd} = \frac12 (g_{ac}g_{bd} -
g_{ad}g_{bc})$.

Due to the anti-symmetry of the indices, any $F \in \Lambda_m M$ will have even rank when expressed as a matrix. Furthermore, as
the dimension of the manifold is four, the rank of any non-zero member of $\Lambda_m M$ is two or four. If the rank of $F$ is
two we say the bivector $F$ is {\it simple}, while if the rank is four $F$ is called {\it non-simple}. If $F$ is simple one may
write $F^{ab} = p^a q^b - q^a p ^b = 2 p^a \wedge q^b$, where $p,q \in T_m M$. By algebraically classifying the simple and
non-simple elements, one is able to identify the possible subalgebras in terms of two simpler subalgebras, $S_m^+ = \{ F \in
\Lambda_m M : F^* = F\}$ and $S_m^- = \{ F \in \Lambda_m M : F^* = -F\}$ \cite{Wang}, where $F^*$ denotes the dual operator.

Noting that $\Lambda_m M = S^+_m \oplus S^-_m$, the authors enumerate all possible subalgebras of $\mathfrak{o}(2,2)$ in
bivector form by examining the subalgebras of $S^+_m$ and $S^-_m$, producing a list of potential subalgebras of
$\mathfrak{o}(2,2)$ \cite{Wang}, with basis vectors taken from the following list of bivectors: \beq & F_1 =
\frac12 (l \wedge n - L \wedge N),~F_2 = \frac12 (l \wedge N),~~F_3 = \frac12 (n \wedge L); & \nonumber \\ &G_1 = \frac12 (l
\wedge n + L \wedge N),~G_2 = \frac12 (l \wedge L),~~G_3 = \frac12(n \wedge N). & \nonumber \eeq \noindent Here, $F_i \in S^+$
and $G_i \in S^-$ for $i\in [1,3]$.
 
With all possible subalgebras of $\mathfrak{o}(2,2)$ identified, we may consider the holonomy group of $(M,g)$, $\Phi$, and
holonomy algebra $\phi$. We will study a particular subalgebra of $\phi$ at each point in the manifold, the infinitesimal holonomy algebra at $m\in M$, $\phi_m'$, arising from contractions of the curvature tensor, $R^a_{~bcd}X^c Y^d,~R^a_{~bcd;e}X^cY^dZ^3$, and so on; where $X,Y,Z,...\in T_m M$. The unique connected group generated by $\phi_m'$ is  a subgroup of the holonomy group $\Phi$ at each point of the manifold, and it is known as the infinitesimal holonomy group $\Phi'_m$ \cite{Hall}. If (M,g) is simply connected and analytic $\Phi'_m$ for each point in M, and $\Phi$ will coincide. 

In general, the bivector representation of $\phi_m'$ as a Lie subalgebra of $\mathfrak{o}(2,2)$ is important as well. This is due to the Ambrose-Singer theorem which states that by computing all of the curvature two-forms from the curvature tensor, and parallel transporting them to a point, $m$, of the manifold, the resulting  collection of two-forms at $m$ span the holonomy algebra. 

The classification by Ghanam and Thompson \cite{GT} requires that the manifold,  metric and connection are analytic to ensure that the Lie algebra of the holonomy group can be computed point-wise from the curvature tensor and its covariant derivatives. In the current work we will impose this condition on the manifold, metric and connection for this reason. 

With the holonomy algebras we may identify any recurrent vectors that are 	admitted by the metric. If there exists a vector $0 \neq {\bf k} \in T_m M$ such that ${\bf k}$ is an eigenvector of each member of $\phi$, then m admits a coordinate neighbourhood $U$, and a nowhere zero vector field ${\bf K}$ on $U$ which agrees with ${\bf K}$ at m, and is such that ${\bf K}$ is recurrent on U \cite{Hall}. 

\begin{subsection}{Example 1} We calculate the infinitesimal holonomy algebra in the null tetrad basis (\ref{Coframe1} - \ref{Coframe2}) with metric functions satisfying (\ref{eqn:RicciFlat1} - \ref{eqn:RicciFlat2}). To do so, we need only contract the Riemann tensor with bivectors constructed from the null tetrad vectors \cite{Nathan}, as the covariant derivatives of the Riemann tensor do not introduce any new elements of the Lie algebra. The matrices with indices lowered are presented, since these have a simpler form: \beq &R_{abcd}
{l_1}^c {n_1}^d =\xi (l_1 \wedge l_2)_{ab}& \\ &R_{abcd} {l_1}^c {l_2}^d = R_{abcd} {l_1}^c {n_2}^d =R_{abcd} {n_1}^c {l_2}^d
=0&\\ &R_{abcd} {n_1}^c {n_2}^d = \xi (l_1 \wedge n_1 - l_2 \wedge n_2 - \zeta ( l_1 \wedge l_2 ) )_{ab}&\\ &R_{abcd} {l_2}^c
{n_2}^d = \xi (l_1 \wedge l_2)_{ab},& \eeq
\noindent where $\xi = \frac{c_1^2 u^2 - 2\alpha c_2 x^2}{2 c_1 u^2 x^2}$ and $\zeta$ is a complicated expression. 

Taking linear combinations of these, we find that our infinitesimal holonomy algebra is spanned by $\left\{ l_1 \wedge n_1 - l_2 \wedge n_2 , l_1 \wedge l_2\right\} $ at each point in the manifold, corresponding to Wang and Hall's subalgebra 2(d) \cite{Wang}. For the holonomy subalgebra $2(d)$, 
$\phi = <F_1, G_2>$, with $|F_1|=-1$ and $|G_2|=0$. 

%
%

From theorem 8.6 in \cite{Hall}, this metric admits two recurrent vectors $\ell$ and $L$ as these are shared eigenvectors of
$F_1$ and $G_2$, with differing eigenvalues. Due to the symmetrization of the two-form representations of the Lie algebra
members in the Riemann tensor, this implies the vectors $\ell$ and $L$ may be seen as eigenvectors with "eigen two-forms"
proportional to the Lie algebra members.


As matrices in the null tetrad basis with the first index up, we may use the Jordan normal form to show  that these two matrices are equivalent to the subcase $A_{10}$ in \cite{GT}. We have found a holonomy algebra isomorphic to $A_{10}$ which, corresponds to a metric containing two recurrent vectors \cite{GT}. A vector $l$ is called recurrent if $\nabla l = l\otimes \omega$ for some one-form $\omega$
\cite{Hall}.

It is worthwhile to consider when the infinitesimal holonomy algebra becomes one-dimensional at all points in the manifold, that is, when $\xi = \frac{c_1^2
u^2 - 2\alpha c_2 x^2}{2 c_1 u^2 x^2} = 0$. This occurs when $c_1$ vanishes and either $a=0$ or $c_2 = 0$. There are three
possible cases where this can happen \begin{enumerate} \item Case 1: $c_1 = \alpha = \beta = a = 0$ 
\item Case 2: $c_1 = \alpha = \beta = c_2 = 0$ 
\item Case 3: $c_1 =c_2 = d = \beta = 0$ 
\end{enumerate}

\noindent According to theorem 4.6 in \cite{GT}, as each of these subcases admit a one-dimensional holonomy
algebra $A_9$, each of these subcases admit two covariantly constant null vectors. This condition implies that these Walker
metrics are Kundt with two null geodesic, expansion-free, shear-free and vorticity-free vectors. For this particular example,
this implies there is a coordinate system where $A_2$ and $B_{01}$ both vanish. Looking at the five examples, it is not clear
that one has found the appropriate coordinate system as only one of $A_2$ or $B_{01}$ vanishes in cases 1,2,4 and 5, while in
case 3 neither function vanishes. This question cannot be answered by comparing holonomy algebras alone, one must
examine these subcases in the context of the equivalence algorithm.

\end{subsection}

\begin{subsection}{Example 2}

With the Riemann components computed, we may contract with frame vectors to determine the infinitesimal holonomy Lie algebra. Despite the differing order of variables, i.e., $\{V,v,U,u\}$ instead of $\{u,v,U,V\}$, we may compare the matrices
arising from the curvature tensor with those in \cite{GT}. We note that the covariant derivatives of the Riemann tensor introduce
no new members of the Lie algebra.

When $B_{10,u} \neq 0$ the Lie algebra is three dimensional, one may use the Jordan canonical form to show that this
is equivalent to $A_{26}$ with $\alpha =0$ in \cite{GT}.
%
\noindent If $B_{10,u} = 0$ there is a bifurcation and the Lie algebra is two dimensional, and hence is equivalent to $A_{17}$
in \cite{GT} 
\noindent Finally, if $B_{10,u} = \Psi = 0$ the Lie algebra is equivalent to the one-dimensional Lie algebra $A_{9}$.

According to Table II in \cite{GT}, when $B_{02}\neq0$, $\Psi\neq0$, and $B_{10_u}\neq0$, we find that the infinitesimal holonomy algebra corresponds to $A_{26}$ at each point in the manifold, and so our metric has a null two-dimensional distribution containing a recurrent
vector field. When $B_{02}\neq0$, $\Psi\neq0$, and $B_{10_u}=0$, we find that the  holonomy algebra corresponds to
$A_{17}$, implying that our metric has a null two-dimensional distribution containing one parallel vector field. We expect that $\partial_V $ will be recurrent because it is an eigenvector of each of the members of the holonomy algebra. Calculating the covariant derivative of $\partial_V$, we find that $\nabla \partial_V = B_{10} \partial_V dU$, as expected. It is also clear that $\partial_V$ becomes parallel when $B_{10} = 0$. If $B_{10,u}$ and $\Psi$ both vanish, the null two-dimensional distribution contains two parallel vector fields.

\end{subsection}

\end{section}

\begin{section}{Holonomy and Recurrent Vectors} Theorem 8.6 in  \cite{Hall} suggests that in order to find the
recurrent vectors, we should calculate the eigenvectors of the elements of the holonomy algebra. In the tetrad basis, we
find that each of the tetrad vectors is an eigenvector of $l_1 \wedge n_1 + l_2 \wedge n_2$, but only $l_1 = \partial_v$ and
$l_2 = \partial_V$ are eigenvectors of $ l_1 \wedge l_2$. The theorem tells us that for each $m\in M$ there exist two recurrent
vector fields on $M$: one having the value $\partial_v$ at $m$ and one having the value $\partial_V$ at $m$. Thus, we look for
two recurrent vectors: one of the form $\boldsymbol\ell_1 = \partial_v+f(u,v,U,V)\partial_V$ and another of the form
$\boldsymbol\ell_2 = h(u,v,U,V) \partial_v + \partial_V$.

We first take the covariant derivative of $\partial_v+f(u,v,U,V)\partial_V$ to find conditions on $f$ that give recurrence for
$\nabla( \ell_1 )$ as: \small

\beq & \alpha\frac{f (d u U^2 + c_1 U) + d u^2 U + u c_2}{c_1 u U} \partial_v \otimes du+ \frac{2f\alpha(u^2 d U + u c_2) +
c_1^2 U + c_1 d u U^2}{2 c_1 u U} \partial_v \otimes dU + & \nonumber \\ & \frac{2 f \alpha (d u^2 U + u c_2) + c_1^2 U + c_1 d
u U^2 + 2 f_u c_1 u U}{2 c_1 u U} \partial_V \otimes du + \frac{f(c_1 U + d u U^2) + c_2 u + d u^2 U + 2 f_U u U}{2 u U}
\partial_V \otimes dU & \nonumber \\ & + f_v \partial_V \otimes dv + f_V \partial_V \otimes dV. & \nonumber \eeq \normalsize

\noindent Thus we must have $f(u,v,U,V)=f(u,U)$ and $f(u,U)$ must satisfy a system of two partial differential equations:
\beq &\alpha\frac{f (d u U^2 + c_1 U) + d u^2 U + u c_2}{c_1 u U} = \frac{1}{f}\frac{2 f \alpha (d u^2 U + u c_2) + c_1^2 U +
c_1 d u U^2 + 2 f_u c_1 u U}{2 c_1 u U} & \nonumber\\ & \frac{2f\alpha(u^2 d U + u c_2) + c_1^2 U + c_1 d u U^2}{2 c_1 u U} =
\frac{1}{f}\frac{f(c_1 U + d u U^2) + c_2 u + d u^2 U + 2 f_U u U}{2 u U}. & \eeq \normalsize

\noindent After some algebra, these become \beq &f^2\alpha\left(\frac{d u U}{c_1} + 1 \right)-\frac{c_1}{2}-\frac{d u U}{2} - u
f_u = 0, & \label{eqn:recrnt0} \\ & f^2\alpha\left(\frac{d u U}{c_1}+\frac{c_2}{c_1}\right)-\frac{c_2}{2}-\frac{d u U}{2} - U
f_U = 0. & \label{eqn:recrnt1} \eeq

Assuming $\alpha\neq0$, these equations have a solution: \beq &f_{k_0}(u,U)=-\sqrt{\frac{c_1}{2\alpha}} \tanh\left(
\sqrt{\frac{\alpha}{2 c_1}}\left( c_1 \ln|u| + c_2 \ln|U| +d u U + k_0 \right) \right) & \nonumber \eeq

\noindent where $k_0$ is an arbitrary constant. This solution gives a one-parameter family of recurrent vectors; if we choose
any $k_1$ and $k_2$ such that $k_1 \neq k_2$, then $f_{k_1}(u,U) \neq f_{k_2}(u,U)$, and so the recurrent vectors
$\partial_v+f_{k_1}(u,U)\partial_V$ and $\partial_v+f_{k_2}(u,U)\partial_V$ are linearly independent.
Repeating the above for $\boldsymbol\ell_2$, we find $h_{k_0'}(u,U) = \frac{2 \alpha}{c_1} f_{k_0'}(u,U)$. However,
choosing $k_0' = k_0+i \pi\sqrt{\frac{c_1}{2 \alpha}} \in\mathbb{C}$ and noting that $\tanh(x+i\pi/2) = 1/\tanh(x)$, we find
that each $h_{k_0'}(u,U)$ corresponds to $h_{k_0}(u, U) = \frac{1}{f_{k_0}(u,U)}$. These solutions are not quite as
different from those for $\boldsymbol\ell_1$ as they initially appear.  If $m \in M$ corresponds to $(u_0, v_0, U_0, V_0)$, then choosing $k_0 = - c_1 \ln|u_0| + c_2 \ln|U_0| + d u_0 U_0$ and $k_0' = k_0$ gives $f_{k_0}(u_0, v_0, U_0, V_0) = h_{k_0'}(u_0,  v_0, U_0, V_0) = 0$ and so our recurrent vectors are such that $\left.\left[\partial_v + f_{k_0} \partial_V\right]\right|_m =
\partial_v$ and $\left. \left[h_{k_0}\partial_v + \partial_V\right]\right|_m = \partial_V$, as Theorem 8.6 \cite{Hall} predicts.

When $\alpha=0$, \eqref{eqn:recrnt0} and \eqref{eqn:recrnt1} have the solution \beq & \widetilde{f_{k_0}}(u,U) = -\frac{ d u U +
c_1 \ln|u| + c_2 \ln|U| + k_0}{2}= \lim_{\alpha\to 0} f_{k_0} & \label{eqn:fKundt}. \eeq

\noindent The corresponding partial differential equations for $h$ have the following  solution \beq & \widetilde{h_{k_0}}(u,U) =
\frac{1}{\widetilde{f_{k_0}}(u,U)} = \lim_{\alpha\to 0} h_{k_0}. & \label{eqn:hKundt} \eeq

\begin{subsection}{Recurrent Vectors in a Kundt Subcase} As a simple example, we show that in the Kundt subcase where $\alpha =
0$, the tetrad vector $l_2 = \partial_V$ is a Kundt vector. We know that $l_2$ is null since it is a null tetrad member.
$\nabla_{l_2} l_2= 0$ implies that $l_2$ is geodesic. It is a simple task to verify that the covariant derivative of $l_2$ is of the form \eqref{KundtVect}

The above recurrent and Kundt vectors were obtained under the assumption (from \eqref{RicciB02eq0}) that $a = \alpha d/c_1$
(including $c_1 = c_2$). Now the only case not examined is $a \neq \alpha d/c_1$ and $c_1 = c_2$ ($c_1 \neq 0$). Resuming from
\eqref{RicciB02eq0} and assuming $c_1 = c_2$ instead of $a = \alpha d/c_1$ results in the same holonomy algebras.
The recurrent vectors are found similarly, with $f(u, U)$ and $h(u, U)$ satisfying the partial differential equations: 
\beq & f^2 \left( a U + \frac{\alpha}{u} \right) - \frac{c_2}{2 u} - \frac{d U}{2} - f_u =0,~~ f^2 \left( a u + \frac{\alpha}{U}
\right) - \frac{c_2}{2 U} - \frac{d u}{2} - f_U = 0, & \nonumber \\ &h^2 \frac{1}{2} \left( d U + \frac{c_2}{u}
\right) - a U - \frac{\alpha}{u} - H_u = 0,~~ h^2 \frac{1}{2} \left( d u + \frac{c_2}{U} \right) - a u - \frac{\alpha}{U} - H_U
= 0. & \nonumber \eeq 

Evidently, when $A_2=0$ ($\Leftrightarrow a = \alpha = 0$), these differential equations reduce to those found previously
and so the recurrent vectors are the same as \eqref{eqn:fKundt} and \eqref{eqn:hKundt}, so we find that we have a Kundt vector
once again.

\end{subsection} \end{section} 

\begin{section}{ The Equivalence Algorithm for a Kundt Neutral Signature Walker Metric}

As in the case of the Lorentzian signature, the equivalence of neutral VSI metrics in general, may  be determined using the Cartan algorithm \cite{mcnutt,Dario}. 
The goal of the equivalence algorithm is the computation of a finite list of invariants arising from the curvature tensor and its 
covariant derivatives which has been normalized by fixing all frame transformations affecting the form of the curvature tensor and 
its covariant derivatives. To begin the equivalence algorithm for this particular class of neutral metrics, we determine the effect 
of the frame transformations on the Riemann tensor. This may be done by computing the effect of each null rotation about
$\ell_1, n_1, \ell_2$ and $n_2$ along with the effect of boosts in the $(\ell_1, n_1)$ and $(\ell_2, n_2)$ planes on the 
spin-coefficients, as these quantities are of vital importance when computing the covariant derivatives of the curvature tensor.

For the Ricci-flat $VSI$ Walker metrics, we may always fix the boost parameters so that two components of the curvature tensor are 
constant. Thus the null rotations are left as potential members of the zeroth order isotropy group. After recording the number of functionally independent invariants
that appear at zeroth order, we proceed to compute the first covariant derivative of the curvature tensor.

From the components of this rank five tensor we may determine the first order Cartan invariants by fixing all frame transformations that affect the form of the first order covariant derivative of the curvature tensor, and identify all new functionally independent and dependent invariants that appear at first order after this process. The algorithm continues each iteration by computing higher order covariant derivatives of the tensor and identifying the isotropy group and functionally independent invariants at each order. The algorithm stops when it reaches the $q$-th iteration for which the dimension of the isotropy group and number of functionally independent invariants does not change from iteration $q-1$ to $q$.

In general it is not known how many iterations are required to compute the entire list of invariants for a particular metric. The theoretical upper-bound introduced by Cartan limits the number of iterations required to classify an arbitrary neutral metric; this upper-bound is determined by the largest isotropy subgroup of the Riemann curvature tensor $\tilde{s}_0$ which must be less than six for spacetimes which are not locally homogeneous: \beq q \leq n + \tilde{s}_0 +1 = 4+5+1 = 10 \nonumber \eeq \noindent In
the case of the Ricci-flat VSI neutral metrics, we may fix the two real-valued boost parameters to set two components of the
curvature tensor to be constant, this reduces the upper-bound from ten to eight. In the Ricci-flat VSI Kundt-Walker metric this
may be reduced further as the isotropy group consists of the two-dimensional null rotations about a particular null vector 
\beq q \leq n + s_0 + 1 = 4+2+1 = 7. \nonumber \eeq 

\noindent This is the standard upper-bound for Lorentzian metrics, however this has only been achieved for a simple subcase of the 
neutral metrics. An effective lowering of the upper-bound for all neutral metrics would require a classification akin to the 
Petrov classification for Lorentzian metrics.

We now describe in detail the equivalence algorithm outline in \cite{CCMMB} for those Ricci-flat Walker metrics with $B_{02}
\neq 0$ and $A_2 = C_2 =0 $ with the following conditions on the remaining metric functions. \beq & B_{10}=f(U), B_{00}= 0,
B_{02}= e^{W(u)}e^{Z(U)}.&\nonumber \eeq \noindent The metric functions $A_0$, $A_1$ and $C_{11}$ now become: \beq A_{0} &=&
\frac18 \frac{-2B_{10} C_{11} - 4 A_1 B_{01} - 2 C_{11,U} + C_{11}^2 }{B_{02}} \nonumber \\ A_{1} &=& \frac12 W_{,u} \nonumber
\\ C_{11} &=& 2 B_{10} + Z_{,U} + G(U). \nonumber \eeq \noindent The non-zero components of the Riemann tensor are: \beq
R_{2323} &=& 2 B_{02} \nonumber \\ R_{2424} &=&-C_{0,U u}+B_{00,u u}-3 v^2 A_{1} B_{02,u}-2 v^2 A_{1,u} B_{02}-B_{10} v
C_{11,u}\nonumber \\ &&+B_{10} v A_{1,U} -A_{1} B_{10} C_{0}+2 v^2 A_{1}^2 B_{02}+B_{10} C_{11} A_{0}\nonumber \\ &&+(1/2) v
C_{11} C_{11,u} -4 A_{0} v B_{02,u}-2 A_{0,u} v B_{02}-v A_{1} B_{01,u}-v A_{1,u} B_{01}\nonumber \\ &&+4 v A_{1} B_{02}
A_{0}+A_{0,U U} +v^2 B_{02,u u}+V B_{10,u u}-2 A_{0} B_{01,u}-A_{0,u} B_{01}\nonumber \\ &&-B_{10} C_{0,u} -A_{1} C_{0,U}
-A_{1,U} C_{0}+v A_{1,U, U}+C_{11} A_{0,U}+C_{11,U} A_{0}\nonumber \\ &&+2 A_{0}^2 B_{02}+(1/2) C_{11,u} C_{0} -v C_{11,U u}+v
B_{01,u u}+A_{1} B_{00,u}+B_{10} A_{0,U}, \nonumber \eeq

\noindent where again we will denote $R_{2424}$ as $\Psi$.

Since we may fix the components of the Riemann curvature tensor to constants by performing boosts in both the $(\ell_1,n_1)$ and
$(\ell_2, n_2)$ null planes, with $z_{1}=A$, $z_{2}=B$ (see Appendix for transformation rules) as boost parameters: \beq z_1^2 =
\frac{z_2^2}{2B_{02}},~~z_2^4 = \frac{2B_{02}}{\Psi} \label{Walkerpp} \eeq \noindent no new functionally independent invariants
appear at zeroth order \cite{mcnutt,Dario}. Thus we must compute the first covariant derivative of the curvature tensor, which
requires knowledge of the spin-coefficients.

The non-vanishing spin coefficients arising from the metric coframe are:
\small
\beq
 \gamma&=&f(U) +\frac{1}{4}(G(U)+Z_{,U})  \nonumber \\
\sigma' &=& \frac12 C_{11} A_{0} - \frac12 v C_{11,u} - \frac12 C_{0,u}+ v A_{1,U}+A_{0,U} - \frac12 A_1 C_0,  \nonumber \\
\kappa' &=& -\frac12 B_{10} v C_{11} - \frac12 B_{10}C_0-2 v^2 A_1 B_{02} - v A_1 B_{01} - 2 A_0 v B_{02} - A_0 B_{01} + V B_{10,u}, \nonumber \\ 
&&+ v^2 B_{02,u} + v B_{01,u}  + v B_{01,u}+B_{00,u} -\frac12 v C_{11,U}-\frac12 C_{0,U} + \frac14 v C_{11}^2 + \frac14 C_{11} C_0, \nonumber\\
\beta' &=& \frac14 W_{,u}, \nonumber \\
\tilde{\beta} &=& - \frac14 W_{,u}, \nonumber \\
\tilde{\gamma} &=& -\frac12 (Z_{,U} + G(U)) ,\nonumber \\
\tilde{\rho}' &=& - f(U)-\frac12 f_{,U} - \frac12 G(u) \nonumber \\
\tilde{\kappa}' &=& 2ve^{W(u)+f(U)}+B_{01} \nonumber \\
\nonumber \eeq
\normalsize 

 We notice that, as we have chosen that $\ell_a n^a = 1$ and $m_a \tilde{m}^a =
-1$, we have the following relationships between spin-coefficients \beq & \epsilon = -\gamma', \alpha = \beta', \beta = \alpha',
\gamma = -\epsilon', \tilde{\epsilon} = -\tilde{\gamma}', \tilde{\alpha} = \tilde{\beta}', \tilde{\beta}=\tilde{\alpha}' ,
\tilde{\gamma} = - \tilde{\epsilon}'. & \label{LawCond} \eeq \noindent Thus, of the thirty-two spin-coefficients we may concentrate on
twenty-four of them instead.

Performing boosts in both the $(\ell_1,n_1)$ and $(\ell_2,n_2)$ null planes and denoting our boosted spin coefficients with a
subscript B (where $\kappa_{B}$ would be the spin coefficient $\kappa$ after boosts in the $(\ell_1, n_1)$ and $(\ell_2,n_2)$
null planes, we find the non-vanishing transformed spin coefficients to be \small \beq
&\alpha'_{B}=-\frac{1}{4}\left(\frac{z_{1,U}}{z_{1}z_{2}}+\frac{z_{2,U}}{z_{2}^{2}}\right),~\alpha_{B}=\frac{1}{2}z_{2}\alpha,~\tilde{\alpha}'_{B}=\frac{1}{2}z_{2}\tilde{\alpha}',~\tilde{\alpha}_{B}=\frac{1}{4}\left(\frac{z_{2,U}}{z_{2}^{2}}-\frac{z_{1,U}}{z_{1}z_{2}}\right)
& \label{Boost1} \\ &
\gamma_{B}=z_{1}\gamma,~\gamma_{B}'=\frac{1}{4}\left(\frac{z_{1,u}}{z_{1}^{2}}+\frac{z_{2,u}}{z_{1}z_{2}}\right),~\tilde{\gamma}_{B}=\tilde{\gamma},~
\tilde{\gamma}_{B}'=\frac{1}{4}\left(\frac{z_{1,u}}{z_{1}^{2}}-\frac{z_{2,u}}{z_{1}z_{2}}\right) & \label{Boosts2}\\
&\sigma_{B}'=\frac{1}{2}z_{1}z_{2}^{2}\sigma',~ \kappa_{B}'=\frac{1}{2}z_{1}^{2}z_{2}\kappa',~
\tilde{\rho}_{B}'=\frac{1}{2}z_{1}\tilde{\rho}',~ \tilde{\kappa}_{B}'=\frac{1}{2}\frac{z_{1}^{2}\tilde{\kappa}'}{z_{2}^{2}}&
\label{Boost4} \eeq

\normalsize

\noindent From the components of this rank-five tensor, we may solve for the following boosted spin-coefficients [10] as first-order
Cartan invariants: \beq & \{ \rho, \tau, \kappa, \sigma, \tilde{\rho}, \tilde{\sigma}, \tilde{\tau}, \tilde{\kappa}, \alpha,
\alpha',\tilde{\alpha},\tilde{\alpha}', \gamma, \gamma', \tilde{\gamma}, \tilde{\gamma}' \} & \nonumber \eeq \noindent Of which,
the following are non-zero: \beq & \{ \alpha, \alpha',\tilde{\alpha},\tilde{\alpha}', \gamma, \gamma', \tilde{\gamma},
\tilde{\gamma}'\}& \nonumber \eeq 

The remaining isotropy at first order consists of null rotations about $\ell_{1}$,
as null rotations about $n_{1}$, $\ell_{2}$, and $n_{2}$ change the number of non-zero components of the Riemann tensor, and
thus do not belong to the first-order isotropy group. Computing null rotations about $\ell_{1}$ and denoting boosted and rotated
spin coefficients with a subscript R, we obtain the following list of non-vanishing transformed spin coefficients (with $z_{3}$,
$z_{4}$ rotation parameters): \beq \begin{array}{c} \alpha_{R}=\alpha_{B}+z_{3}\gamma_{B}'\\
\alpha_{R}'=\alpha_{B}'-z_{4}\gamma_{B}'\\ \tilde{\alpha}_{R}=\tilde{\alpha_{B}}-z_{4}\tilde{\gamma_{B}}'\\ \tilde{\alpha}_{R}'=
\tilde{\alpha_{B}}'+z_{3}\tilde{\gamma_{B}}' \end{array} \nonumber ~
\begin{array}{c}\gamma_{R}=z_{3}z_{4}\gamma_{B}'-z_{3}\alpha_{B}'+z_{4}\beta_{B}'+\gamma_{B}\\ \gamma_{R}'=\gamma_{B}'\\
\tilde{\gamma}_{R}=z_{3}z_{4}\tilde{\gamma}_{B}'-z_{3}\tilde{\beta}_{B}'+z_{4}\tilde{\alpha}_{B}'+\tilde{\gamma}_{B}\\
\tilde{\gamma}_{R}' =\tilde{\gamma}_{B}'\end{array} \eeq\beq & \tau_R = 0,~~\tilde{\tau}_R = 0,~~\rho_R = 0,~~\tilde{\rho}_R =
0,~~\sigma_R = 0,~~\tilde{\sigma}_R = 0,~~ \kappa_R = 0,~~\tilde{\kappa}_R = 0 \nonumber \eeq

As $\gamma_{B}'$ and $\tilde{\gamma_{B}}'$ are unaffected by the null rotation, they are invariant under such a transformation;
that is, $\gamma_{B}'=\gamma_{R}'$ and $\tilde{\gamma_{B}}'=\tilde{\gamma_{R}}'$. Furthermore, their vanishing or non-vanishing
affects the transformation rules for the remaining first order invariants, and hence indicates possible subcases. 

In the present work, we examine a simple subcase in order to present a complete application of the equivalence algorithm. We will assume that the following spin-coefficients are equal to zero: 
\beq \{ \alpha_B, \alpha_B' , \tilde{\alpha}_B , \tilde{\alpha}_B' , \gamma_B' ,
\tilde{\gamma}_B' \} .\nonumber \eeq 

The vanishing of these spin-coefficients produce the following conditions on the metric functions: \beq z_{1,u}=z_{1,U} = z_{2,u}=z_{2,U} = 0. \nonumber \eeq
\noindent The components of the Riemann tensor \eqref{Walkerpp} are constant, implying that $Z_{,U} = W_{,u} = 0$, and so
the metric function $A_1$ vanishes. The constancy of the curvature component $\Psi$ requires that that $C_0$ satisfies a
complicated partial differential equation.

Simplifying the above expressions for our boosted and rotated spin coefficients, we obtain two non-zero first order invariants:
\beq \{ \gamma_B, \tilde{\gamma}_B ' \}.\nonumber \eeq However, since
$\beta_{B}'=\alpha_{B}'=\tilde{\alpha_{B}}'=\tilde{\beta_{B}}'=0$, our new first order invariants are unchanged under null
rotations, and so we cannot fix all of our isotropy after first order. The dimension of the isotropy group after first order
$dim(I)=2$, and we must proceed to second order.

Taking the second order covariant derivative of the Riemann tensor, we may simplify the components of this rank six tensor
to produce the following set of second order curvature invariants \beq \{ \gamma_B, \tilde{\gamma}_B,
D\gamma_{B},\delta\gamma_{B}, D\tilde{\gamma_{B}}, \delta\tilde{\gamma_{B}}, \Delta\tilde{\gamma_{B}}, \Delta\gamma_{B},
D'\gamma_{B}, D'\tilde{\gamma_{B}} \} \nonumber \eeq \noindent Since we have that $\gamma_{B}$ and $\tilde{\gamma}_{B}$ are
functions of U alone, all of their derivatives taken with respect to u, v, and V vanish, the frame derivatives
simplify to
$D\gamma_{B}=\delta\gamma_{B}=D\tilde{\gamma_{B}}=\delta\tilde{\gamma_{B}}=\Delta\gamma_{B}=\Delta\tilde{\gamma_{B}}=0$, and
$D'\gamma_{B}=\partial_{U}\gamma_{B}, D'\tilde{\gamma_{B}}=\partial_{U}\tilde{\gamma_{B}}$. 
The scalars $\gamma_{B}$ and
$\tilde{\gamma_{B}}$ are invariants, they are expressions not involving $z_{3}$ or $z_{4}$; that is, \begin{align*}
D'\gamma_{B}&=D'\left[\gamma_{B}\right]=\partial_{U}\left[\gamma_{B}\right],\\
D'\tilde{\gamma_{B}}&=D'\left[\tilde{\gamma_{B}}\right]=\partial_{U}\left[\tilde{\gamma_{B}}\right]. \end{align*} We
cannot manipulate these equations to find conditions on $z_{3}$ and $z_{4}$, and so the
dimension of the isotropy group after second order remains $2$. The algorithm terminates at the second iteration, as
$t_1 = t_2 = 1$ and $dim H_1 = dim H_2 = 2$. The resulting list of invariants up to second order allows one to completely
classify these spaces, and no further iterations of the algorithm will yield no new information.

%
\end{section} 

\section{Discussion}

We have investigated the mathematical properties of a class of four-dimensional neutral signature metrics, with vanishing scalar
curvature invariants (VSI). We examined a collection of metrics which satisfy the $VSI$-property and are distinct from the Kundt class. To discuss the difference in the neutral Ricci-flat Walker metrics with vanishing scalar curvature invariants, we compared two analytic metrics with different 2-dimensional holonomy algebras: one which is generally Walker but not Kundt, and a second that is always Kundt.

By giving conditions for the existence of a null geodesic, expansion-free, shear-free, and vorticity-free vector for Walker
metrics we were able to compare the two examples. Then, using the Lie algebra classification provided in \cite{GT} and
\cite{Wang}, we explicitly identified the geometrically special vectors that arise from the holonomy algebra in
each example. This classification is well-suited for determining the existence of invariant null distributions, recurrent
vectors and covariantly constant null vectors; however it is not fine enough to determine the equivalence of metrics. As an
example it is clear that the two metrics are inequivalent as they have distinct two-dimensional Lie algebras, yet both metrics
contain a subcase for which these Lie algebras become one-dimensional. This is notable as all one-dimensional Lie algebras are
equivalent to the Lie algebra for a metric admitting two covariantly constant null vectors, implying that this metric is doubly
Kundt.

A natural question is to ask when, if at all, are the subcases of the two metrics equivalent. This question can only be resolved
by implementing the equivalence algorithm for neutral metrics, which is a non-trivial task. We have provided a simple example of
the equivalence algorithm applied to a subcase of the Kundt-Walker metric, which parallels the plane-wave spacetimes in the
Lorentzian case. We have shown that neutral signature "plane waves" require the same number of covariant derivatives as their
Lorentzian counterparts. It is unknown whether this holds for neutral-signature metrics in general, due to the difference in the
group of frame transformations, it is possible that the neutral signature metric require a higher number of covariant
derivatives to complete the equivalence algorithm. In the context of the Ricci-flat Walker metrics this is a particularly
relevant question as one cannot simply compare scalar curvature invariants to determine equivalence. \cite{mcnutt}

\section{Appendix: Transformation Rules for Spin Coefficients} Consider the boosts in the two null planes, given by the
following transfomation, \beq \{ n_a, \ell_a, \tilde{m}_a, m_a \} \to \{ A n_a, A^{-1} \ell_a, B \tilde{m}_a, B^{-1} m_a\}
\label{Boosts} \eeq

\noindent the spin-coefficients transform as: 

\small

\beq & \kappa_B = \frac{\kappa}{A^2 B},~ \rho_B = \frac{\rho}{A},~ \sigma_B = \frac{\sigma}{AB^2},~ \tau_B =\frac{\tau}{B},~&
\nonumber \\ &\tau'_B = B \tau',~ \sigma'_B = A B^2 \sigma',~ \rho'_B = A \rho',~ \kappa'_B = A^2 B \kappa',~ & \nonumber \\
&&\nonumber \\ &\gamma'_B = \frac12 \left[ \frac{D(A) + A \gamma' + A \gammat'}{A^2} + \frac{D(B) + B\gamma' - B \gammat'}{AB}
\right],~ & \nonumber \\ &\beta'_B = -\frac12 \left[ \frac{B( \Delta (A) - A \alphat' - A \beta')}{A} + \Delta (B) + B \alphat'
- A \beta' \right],~ & \nonumber \\ &\alpha'_B = -\frac12 \left[ \frac{ \delta (A) - A \alpha' - A \betat'}{AB} - \frac{-\delta
(B) + B\alpha' - B \betat'}{B^2} \right],~ & \nonumber \\ &\epsilon'_B = \frac12 \left[D'(A) + A \epsilon' + A \epsilont' +
\frac{A(D'(B)+B\epsilon' - B \epsilont')}{B} \right],~ & \nonumber \eeq

\beq & \kappat_B = \frac{B\kappat }{A^2},~ \sigmat_B = \frac{B^2\sigmat }{A},~ \rhot_B =\frac{\rhot}{A} \taut_B = B \taut,~&
\nonumber \\ & \taut'_B = \frac{\taut'}{B},~ \rhot'_B =A \rhot',~ \sigmat'_B = \frac{A \sigmat' }{B^2},~ \kappat'_B = \frac{A^2
\kappat' }{B},~& \nonumber \\ &&\nonumber \\ & \gammat'_B = \frac12 \left[ \frac{D(A)+A\gamma' + A \gammat'}{A^2} - \frac{D(B) +
B \gamma' - B \gammat'}{AB} \right],~ & \nonumber\\ & \alphat'_B = -\frac12 \left[ \frac{B(\Delta (A)-A\alphat' - A
\beta')}{A^2} - \Delta (B) - B \alphat' + B \beta' \right],~ & \nonumber \\ & \betat'_B = -\frac12 \left[ \frac{\delta
(A)-A\alpha' - A \betat'}{AB} + \frac{ - \delta (B) + B \alpha' - B \betat'}{B^2} \right],~ & \nonumber \\ & \epsilont'_B =
\frac12 \left[ D'(A)+A\epsilon' + A \epsilont' - \frac{ A(D'(B) + B \epsilon' - B \epsilont')}{B^2} \right],~ & \nonumber \eeq

\normalsize

To produce a rotation about the null vector $\ell^a$ we make the transformation:
\beq \{n_a, \ell_a, \tilde{m}_a, m_a \} \to \{ n_a + \tilde{\mu} \tilde{m}_a - \mu m_a - \mu \tilde{\mu} \ell_a, \ell_a,
\tilde{m}_a -\mu \ell_a, m_a + \tilde{\mu} \ell_a \} \label{lRot} \eeq

\noindent while the spin-coefficients transform as \small \beq & \kappa_R = \kappa,~ \rho_R = \rho-\mu \kappa,~ \sigma_R =
\sigma+\tilde{\mu} \kappa,~ \tau_R =\tau + \tilde{\mu} \rho - \mu \sigma - \mu\tilde{\mu} \kappa ,~& \nonumber \\ 
&\tau'_R = \tau' + D(\mu) - 2 \mu \gamma' + \mu^2 \kappa ,~ & \nonumber \eeq \
\beq &\sigma'_R = \sigma'- \Delta(\mu) + \mu D(\mu) - 2\mu \beta' - 2 \mu^2 \gamma' - \mu^2 \rho + \mu^3 \kappa + \mu \tau',~ & \nonumber \\ &\rho'_R = \rho'- \delta(\mu) - D(\mu) \tilde{\mu} - 2 \mu \alpha' + 2 \mu \tilde{\mu} \gamma' - \mu^2 \sigma - \mu^2 \tilde{\mu} \kappa - \tilde{\mu} \tau',~ & \nonumber \\ &\kappa'_R = \kappa' + \Delta(\mu) \tilde{\mu} - \mu \delta(\mu) + D'(\mu) - \mu \tilde{\mu} D(\mu) + 2 \mu
\tilde{\mu} \beta' - 2 \mu^2 \alpha' - 2 \mu \epsilon' & \nonumber \\ & + 2 \mu^2 \tilde{\mu} \gamma' + \mu^2 \tilde{\mu} \rho -
\mu^3 \sigma + \mu^2 \tau - \mu^3 \tilde{\mu} \kappa - \mu \sigma' + \mu \rho' - \mu \tilde{\mu} \tau', & \nonumber \eeq 
\beq&\gamma'_R = \gamma' - \mu \kappa, & \nonumber \\ &\beta'_R = \beta' + \mu \rho - \mu^2 \kappa + \mu \gamma', &
\nonumber \\ &\alpha'_R =\alpha' + \mu \sigma + \mu \tilde{\mu} \kappa - \tilde{\mu} \gamma', & \nonumber \\ &\epsilon'_R =
\epsilon' - \mu \tilde{\mu} \gamma' + \mu \alpha' - \tilde{\mu} \beta' + \mu^2 \sigma + \mu^2 \tilde{\mu} \kappa - \mu
\tilde{\mu} \rho - \mu \tau, & \nonumber \\ & \kappat_R =\kappat ,~\sigmat_R = \sigmat-\mu \kappat,~\rhot_R =\rhot + \tilde{\mu}
\kappat,~\taut_R = \taut+ \tilde{\mu} \sigmat - \mu \rhot - \mu \tilde{\mu} \kappat,& \nonumber \eeq 
\beq  & \taut'_R = \taut' - D(\tilde{\mu}) + 2 \tilde{\mu} \gammat' + \tilde{\mu}^2 \kappat,& \nonumber \\ & \rhot'_R =\rhot' + \Delta(\tilde{\mu}) - \mu D(\tilde{\mu}) + 2 \tilde{\mu} \alphat' + 2 \mu \tilde{\mu} \gammat' - \tilde{\mu}^2 \sigmat +\mu \tilde{\mu}^2 \kappat + \mu \taut',& \nonumber \\ & \sigmat'_R = \sigmat' + \delta(\tilde{\mu}) + \tilde{\mu} D(\tilde{\mu}) + 2 \tilde{\mu} \betat' - 2 \tilde{\mu}^2 \gammat' - \tilde{\mu}^2 \rhot - \tilde{\mu}^3 \kappat -\tilde{\mu} \taut',& \nonumber \\ & \kappat'_R = \kappa' - \Delta(\tilde{\mu}) \tilde{\mu} + \mu \delta(\tilde{\mu}) - D'(\tilde{\mu}) - \mu \tilde{\mu} D(\tilde{\mu}) - 2 \tilde{\mu}^2 \alphat' + 2 \mu \tilde{\mu} \betat' + 2 \tilde{\mu} \epsilont' & \nonumber \\ &- 2 \mu \tilde{\mu}^2 \gammat' + \tilde{\mu}^3 \sigmat - \mu \tilde{\mu}^2 \rhot +\tilde{\mu}^2 \taut - \mu \tilde{\mu}^3 \kappat - \tilde{\mu} \rhot' + \mu \sigmat' - \mu \tilde{\mu} \taut' ,& \nonumber \eeq 
\beq & \gammat'_R = \gammat' + \tilde{\mu} \kappat, & \nonumber \\ & \alphat'_R = \alphat' - \tilde{\mu} \sigmat + \mu \tilde{\mu} \kappat + \mu \gammat', & \nonumber \\ & \betat'_R = \betat' - \tilde{\mu} \rhot - \tilde{\mu}^2 \kappat - \tilde{\mu} \gammat', & \nonumber \\ & \epsilont'_R = \epsilont' - \mu \tilde{\mu} \gammat' + \mu \betat' - \tilde{\mu} \alphat' + \tilde{\mu}^2 \sigmat - \mu \tilde{\mu} \rhot - \mu \tilde{\mu}^2 \kappat + \tilde{\mu} \taut, & \nonumber \eeq

\normalsize \noindent To generate a rotation about $n_a$ we may apply the prime operation \cite{Law} to the above
spin-coefficients. Notice that $\mu' = -\mu$ and $\tilde{\mu}' = -\tilde{\mu}$, this is reflected in the resulting frame
transformation on $M$

\beq \{\ell_a, n_a, -m_a, -\tilde{m}_a \} \to \{\ell_a - \tilde{\mu} m_a + \mu \tilde{m}_a - \mu \tilde{\mu} n_a, n_a, -m_a -
\mu n_a, -\tilde{m}_a + \tilde{\mu} n_a \} \nonumber \eeq

\noindent to determine the effect of a null rotation about $n$ on the spin-coefficients merely prime the above quantities.
There are twenty-four discrete transformations that will be important, although it is best seen on the level of vectors on $M$
as the interchange of the order of the four null vectors. As an example, consider the following transformation which is relevant
for the subcase of Ricci-flat VSI Walker metrics we have been studying:

\small

\beq \ell^{a \times} = m^a,~~m^{a \times} = \ell^a,~~n^{a \times} = -\tilde{m}^a,~~\tilde{m}^{a \times} = - n^a. \nonumber \eeq
\noindent Noting that the square of this transformation is identity, we may summarize the effect on the spin coefficients as:
\beq & \epsilon^{\times} = \alpha',~ \alpha^{\times} = \epsilon',~~\beta^{\times} = - \gamma',~~\gamma^{\times} = - \beta, &
\nonumber \\ &\kappa^{\times} = - \sigma,~~\rho^{\times} = \tau,~~\tau^{'~\times} = \rho',~~\sigma^{'~ \times} = -\kappa', &
\nonumber \\ &\epsilont^{\times} = -\betat',~~\betat^{\times} = - \epsilont',~~\alphat^{\times} = \gammat',~~\gammat^{\times} =
\alphat' & \nonumber \\ &\kappat^{\times}= \sigmat',~~\sigmat^{\times} = \kappat',~~\rhot^{\times} = -\taut',~~\taut^{\times} =
-\rhot'. & \nonumber \eeq

\normalsize \noindent Although the priming operation leaves the formula unchanged for this example, this may not be the case
with other re-orderings of the coframe.

\section*{Acknowledgements} The work was supported, in part, by NSERC of Canada (AC) and by the David and Faye Sobey Foundation
(NM). DM would like to thank Peter Law for his concise and helpful remarks. \appendix

\end{document}